\documentclass[12pt]{article}
\usepackage{amsmath,amssymb}
\usepackage{amsthm}
\usepackage{enumerate}
\usepackage{mathtools}
\usepackage{tikz}
\usepackage{authblk}
\usepackage{float}
\makeatletter
\def\widebreve#1{\mathop{\vbox{\m@th\ialign{##\crcr\noalign{\kern3\p@}%
      \brevefill\crcr\noalign{\kern3\p@\nointerlineskip}%
      $\hfil\displaystyle{#1}\hfil$\crcr}}}\limits}
\def\brevefill{$\m@th \setbox\z@\hbox{$\braceld$}%
  \bracelu\leaders\vrule \@height\ht\z@ \@depth\z@\hfill\braceru$}
\makeatletter
\newtheorem{theorem}{Theorem}[section]
\newtheorem*{theorem*}{Theorem}
\newtheorem{proposition}[theorem]{Proposition}
\newtheorem{corollary}[theorem]{Corollary}
\newtheorem{lemma}[theorem]{Lemma}

\newtheorem{prop}[theorem]{Proposition}

\theoremstyle{definition}

\newtheorem*{remark*}{Remark}

\newtheorem*{observation*}{}
\makeatletter
\@addtoreset{equation}{section}
\makeatother

\newcommand{\noin}{\noindent}

\providecommand{\AMS}{$\mathcal{A}$\kern-.1667em%
\lower.25em\hbox{$\mathcal{M}$}\kern-.125em$\mathcal{S}$}

\begin{document}

\date{}

\renewcommand\Authands{, }
\title{Conclusive Posets and Torsion in First Homology}
\author[1]{Bekir
DANIŞ}
\author[2]{İsmail Alperen ÖĞÜT\footnote{Corresponding author}}

\affil[1]{Department of Computer Engineering, Aydın Adnan Menderes University, Aydın, Turkey \protect\\ \texttt{bekir.danis@adu.edu.tr}}
\affil[2]{Department of Electronics and Automation\\
 Ankara University, Vocational School of First Organized Industrial Zone, Ankara, Turkey  \protect\\ \texttt{iaogut@ankara.edu.tr}}

\maketitle

\small

\begin{abstract}
\noin We show that a finite poset is conclusive in the sense of Cigler, Jerman and Wojciechowski if and only if the first integral homology group of its order complex is either zero or has positive free rank. Building on the poset splitting  machinery developed by Cianci and Ottina, we show that posets on $13$ elements with four maximal and four minimal elements do not admit torsion in the first homology, which was the only remaining configuration whose torsion behavior was unknown. Consequently, we obtain that the only  inconclusive posets on $13$ elements are the two dual minimal finite models of the real projective plane.  We also provide a rank criterion for the vanishing of the first cohomology over a field.
\smallskip\\
\noin 2020 {\it Mathematics Subject Classification:}
Primary 06A11, Secondary 06A07.

\smallskip
\noin {\it Keywords:} conclusive poset, derivation, minimal finite model, real projective plane, torsion. 
\end{abstract}
\section{Introduction}
Let $P$ be a poset, $G$ an abelian group. A function  on the comparable pairs of $P$ is called transitive if $f(x,z) = f(x,y) + f(y,z)$ whenever $x\leq y\leq z$. It is called potential if there exists a function $\varphi : P\to G$ such that $f(x,y) = \varphi(y)-\varphi(x)$ for every $x\leq y$. 
A finite poset is called \textbf{soluble} if for every $G$, every transitive function is potential. The poset is called \textbf{defective} if it has a non-potential transitive function for every nontrivial $G$. A poset is called \textbf{conclusive} if it is either soluble or defective. Cigler, Jerman and Wojciechowski conjectured that every finite poset is conclusive \cite[Conjecture $3.1$]{CJW22}, and verified the conjecture for posets with at most $10$ elements (\cite[Corollary $6.7$]{CJW22}).

The conjecture can be reformulated in terms of the derivations of the incidence algebra of the poset. Let $k$ be a unital commutative ring. A transitive function $f$ induces a derivation of the incidence algebra, and every $k$-linear derivation is the sum of an inner derivation and one induced by a transitive function. The derivation induced by $f$ is inner precisely when $f$ is potential. Consequently, there is a canonical isomorphism between the quotient of the $k$-module of derivations by the inner derivations and the first Hochschild cohomology group of the incidence algebra:
$$
\dfrac{Der(I(P,k)) }{Inn(I(P,k)) } \cong HH^1(I(P,k), I(P,k))
$$
where $I(P,k)$ denotes the incidence algebra of $P$ over $k$.

In Section $2$, using the above identification, we will show that a finite poset $P$ is inconclusive precisely when the first integral homology group $H_1(\Delta(P);\mathbb{Z})$ of its order complex is a nonzero finite abelian group. This enables us to establish a connection with the minimal finite models of topological spaces. In particular, the first integral homology group of the minimal finite model of the real projective plane $\mathbb{RP}^2$ is isomorphic to $\mathbb{Z}/2\mathbb{Z}$. Therefore, it constitutes a counterexample to the conjecture above. 

The counterexample is of minimal order. In \cite{A14}, Adamaszek proved by computer verification that the smallest poset whose homology contains torsion has $13$ elements, hence the bound on conclusive posets obtained in \cite{CJW22} can be extended to $12$. In \cite{CO18}, by developing a method called poset splitting, Cianci and Ottina obtained a more conceptual proof of this fact. They also classified the $13$-element finite models of $\mathbb{RP}^2$, showing that there are exactly two such posets, dual to each other, up to homeomorphism.  Their classification is focused on separating the posets by the number of minimal, non-extremal and  maximal elements. In the profile $(3,6,4)$, their argument identifies the two dual minimal finite models of $\mathbb{RP}^2$. It rules out torsion in the profile $(4,4,5)$. For the remaining profile $(4,5,4)$, it shows that the nonextremal elements must form an antichain, but the final exclusion of a model of $\mathbb{RP}^2$ uses the Euler characteristic and does not rule out torsion in first homology. This case turns out to be much more complicated and contains an abundance of possibilities for which new methods have to be developed. The purpose of the present paper is to settle this remaining torsion question.

\begin{theorem}\label{thm1.1}
Let $P$ be a poset on $13$ elements. Then $H_1(\Delta(P);\mathbb{Z})$ contains nonzero torsion if and only if $P$ is isomorphic to one of the two dual minimal finite models of $\mathbb{RP}^2$.
\end{theorem}

Although the above results are very precise, they are developed for a particular case. For the general situation, we also obtain two supplementary results regarding outer derivations. Over a field, we give a criterion for the vanishing of the first cohomology in terms of the rank of a particular matrix. We further strengthen a classical result implying the existence of a crown in a poset that admits outer derivations.

The paper is organized as follows. Section 2 characterizes conclusiveness in terms of first integral homology and establishes a graph theoretic criterion for torsion-freeness. In Section 3, we develop a combinatorial framework which reduces the remaining homological problem to a collection of distinct configurations and eliminate many of these cases. Section 4 deals with the remaining resistant configurations and completes the proof of Theorem 1.1. Section 5 gives two general criteria concerning outer derivations.
\section{Conclusiveness and First Homology}
Throughout this section, let $P$ be a finite poset, $k$ a unital commutative ring, and $I(P,k)$ the incidence algebra of $P$ over $k$. We will establish a homological characterization for conclusiveness. We will also obtain a sufficient combinatorial condition for torsion-freeness, which will be one of the main tools for the later sections.

Let $\Delta(P)$ denote the order complex of $P$. Whenever $k$ is a unital commutative ring, we also regard it as its underlying additive group and obtain the $k$-valued transitive and potential functions accordingly. From the discussion in the introduction, the solubility of $P$ is equivalent to the vanishing of the first Hochschild cohomology group $HH^1(I(P,k),I(P,k))$ for every $k$, whereas the defectivity of $P$ is equivalent to $HH^1(I(P,k),I(P,k))\not = 0$ for every nontrivial commutative ring $k$. Moreover, by the identification given in \cite{GS83} between Hochschild and simplicial cohomologies, we also have $HH^1(I(P,k),I(P,k))\cong H^1(\Delta(P);k)$. This point of view enables the utilization of topological machinery.  More generally, for every abelian group $G$, the group of $G$-valued transitive functions modulo potential functions is naturally isomorphic to $H^1(\Delta(P);G)$. 

For the integral homology, we omit the coefficient ring $\mathbb{Z}$ and simply write $H_1(\Delta(P))$ and extend this notation to other integral homologies. The following provides a complete characterization for the notion of conclusiveness. 
\begin{theorem}
Let $P$ be a finite poset and $\Delta(P)$ its order complex. Then \(P\) is conclusive if and only if \(H_1(\Delta(P))\) is either zero or has a positive free rank.
\end{theorem}
\begin{proof}
As stated above, we are to investigate torsion in $H^1(\Delta(P);G)$. By the Universal Coefficient Theorem
$$
H^1(\Delta(P);G)\cong\mathrm{Hom}_\mathbb{Z}(H_1(\Delta(P)),G).
$$
Write $H_1(\Delta(P)) \cong \mathbb{Z}^r\oplus T$ where $T$ is the finite torsion part. Hence
$$
H^1(\Delta(P); G)\cong \mathrm{Hom}_\mathbb{Z}(\mathbb{Z}^r\oplus T,G)\cong G^r\oplus\mathrm{Hom}_\mathbb{Z}(T,G).
$$
For every nontrivial choice of $G$, if $r>0$, then $H^1(\Delta(P);G)\not = 0$ and $P$ is defective, and if $r=0, T=0$, then $H^1(\Delta(P);G)=0$ and $P$ is soluble, so $P$ is conclusive in each of these cases.
Now suppose that $r=0, T\not=0$, so $H_1(\Delta(P))\cong T.$ If $G=\mathbb{Z}$, we have $$H^1(\Delta(P);G) = \mathrm{Hom}_\mathbb{Z}(T,\mathbb{Z}) = 0.$$
On the other hand, choose a prime $p$ dividing $|T|$. Then $\mathrm{Hom}_\mathbb{Z}(T,\mathbb{Z}/p\mathbb{Z})\not=0$, hence $H^1(\Delta(P);\mathbb{Z}/p\mathbb{Z})\not = 0$ so $P$ is not conclusive if $H_1(\Delta(P))$ is a nonzero torsion group.
\end{proof}
It is known that the minimal finite model of $\mathbb{RP}^2$ is a poset on $13$ elements for which $H_1(\Delta(P))\cong\mathbb{Z}/2\mathbb{Z}$. Thus it constitutes a counterexample to Conjecture $3.1$ in \cite{CJW22}.
\begin{figure}[H] 
	\centering
	\begin{tikzpicture}[scale=0.8,
		thick,
		dot/.style={circle, fill=black, inner sep=0pt, minimum size=4pt}
		]
		
		\node[dot, label=below:$n_1$] (n1) at (-4,0) {};
		\node[dot, label=below:$n_2$] (n2) at (0,0) {};
		\node[dot, label=below:$n_3$] (n3) at (4,0) {};
		
		
		\node[dot, label=left:$a_1$] (a1) at (-5, 3) {};
		\node[dot, label=left:$a_2$] (a2) at (-3, 3) {};

		\node[dot, label=left:$a_3$] (a3) at (-1, 3) {};
		\node[dot, label=right:$a_4$] (a4) at (1, 3) {};
		
		\node[dot, label=right:$a_5$] (a5) at (3, 3) {};
		\node[dot, label=right:$a_6$] (a6) at (5, 3) {};
		
		\node[dot, label=above:$m_1$] (m1) at (-3, 6) {};
		\node[dot, label=above:$m_2$] (m2) at (-1, 6) {};
		\node[dot, label=above:$m_3$] (m3) at (1, 6) {};
		\node[dot, label=above:$m_4$] (m4) at (3, 6) {};
		
		\draw (n1) -- (a1); \draw (n2) -- (a2);
		\draw (n1) -- (a2); \draw (n2) -- (a3);
		
		\draw (n2) -- (a4); \draw (n3) -- (a1);
		\draw (n2) -- (a6); \draw (n3) -- (a4);
		
		\draw (n1) -- (a3); \draw (n3) -- (a5);
		\draw (n1) -- (a5); \draw (n3) -- (a6);
		
		
		\draw (a1) -- (m1); \draw (a3) -- (m2); \draw (a5) -- (m3);

		\draw (a1) -- (m2); \draw (a4) -- (m2); \draw (a6) -- (m4);
		
		\draw (a2) -- (m1); \draw (a3) -- (m4); \draw (a6) -- (m1);

		\draw (a2) -- (m3); \draw (a4) -- (m3); \draw (a5) -- (m4);
		
	\end{tikzpicture}
	\caption{The minimal finite model of $\mathbb{RP}^2$ (13 elements).}
\end{figure}
To stay in the spirit of the conjecture, we provide the following outer derivation of the incidence algebra of the above poset over $\mathbb{Z}/2\mathbb{Z}$. To verify, one can check that it has nonzero circulation (see the paragraph above Theorem \ref{thm5.3}) over the cycle $a_1-m_1-a_6-m_4-a_3-m_2-a_1$ and that it satisfies the consistency relations. 
\begin{table}[H]
	\centering
	\renewcommand{\arraystretch}{1.2} 
	\setlength{\tabcolsep}{8pt}      
	
	\[
	\begin{array}{|c|cccc|ccc|}
		\hline
		a_i & m_1 & m_2 & m_3 & m_4 & n_1 & n_2 & n_3 \\
		\hline
		a_1 & 1 & 0 & \cdot & \cdot & 1 & \cdot & 1 \\
		a_2 & 0 & \cdot & 0 & \cdot & 0 & 0 & \cdot \\
		a_3 & \cdot & 0 & \cdot & 0 & 1 & 0 & \cdot \\
		a_4 & \cdot & 0 & 0 & \cdot & \cdot & 0 & 1 \\
		a_5 & \cdot & \cdot & 0 & 1 & 0 & \cdot & 1 \\
		a_6 & 0 & \cdot & \cdot & 0 & \cdot & 0 & 0 \\
		\hline
	\end{array}
	\]
	
	\caption{Table of values for the transitive function inducing the derivation. Entries represent $f(a, m)$ or $f(n, a)$ in $\mathbb{Z}/2\mathbb{Z}$. The dot $\cdot$ indicates that elements are incomparable.}
\end{table}
For any finite set $S$, let $\lbrace e_s: s\in S\rbrace $ denote the basis of the free abelian group $\mathbb{Z}^S$. Also let $\varepsilon_S : \mathbb{Z}^S\to \mathbb{Z}$ denote the augmentation map. Set $A_S=\mathrm{ker}(\varepsilon_S)$ and for $T\subseteq S$, regard $A_T$ as a subgroup of $A_S$.
Having obtained a characterization in terms of torsion, let us now construct a sufficient combinatorial condition for torsion-freeness.
\begin{lemma}\label{lem2.2}
	Let $U$ be a finite set, $\mathcal{J}$ a family of subsets of $U$. Then $A_U/\sum_{J\in\mathcal{J}}A_J$ is free.
	\end{lemma}
	\begin{proof}
		Let $\sim$ be the equivalence relation on $U$ generated by the relations $u\sim v$ whenever $u,v\in J$ for some $J\in\mathcal{J}$. Let $U_1,U_2,\dots, U_k$ denote the equivalence classes of $\sim$. Define 
		$$\varphi:A_U\to A_{\lbrace 1,\dots,k\rbrace}$$
		 by 
		$$\varphi(\sum_{u\in U}a_ue_u)=(\sum_{u\in U_1}a_u,\sum_{u\in U_2}a_u,\dots, \sum_{u\in U_k}a_u).$$ 
		Since $\sum_{u\in U}a_u=0$, the image belongs to the given codomain. It is clear that $\varphi$ is surjective.
		For every $J\in\mathcal{J}$, the group $A_J$ is generated by the differences $e_u-e_v$ with $u,v\in J$. Since each equivalence class $U_i$ is the connected component of the graph where the edges are between elements which lie in a common member of $\mathcal{J}$, every difference $e_u-e_v$ with $u,v\in U_i$ is a sum of differences belonging to groups $A_J$. It follows that
		$$\sum_{J\in\mathcal{J}} A_J = A_{U_1}\oplus A_{U_2}\oplus\dots\oplus A_{U_k}=\mathrm{ker}(\varphi).$$ 
		Hence $	A_U/\sum_{J\in\mathcal{J}}A_J\cong A_{\lbrace 1,\dots ,k\rbrace}$ which is free.
	\end{proof}
	\begin{lemma}\label{lem2.3}
	Let $G$ be a simplicial complex, $K_1,K_2,\dots, K_m$ connected subcomplexes of $G$. For each $r$, let $C_r=v_r\ast K_r$ be the cone on $K_r$, where $v_r$ are distinct vertices that do not belong to $G$. Set $Y = G\cup\cup_{r=1}^m C_r$. For each $r$, let $\iota_r:K_r\hookrightarrow G$ denote the inclusion map. Then
$$H_1(Y)\cong H_1(G)/\sum_{r=1}^m(\iota_r)_\ast (H_1(K_r)).$$
	\end{lemma}
	\begin{proof}
		Let $Y_0 = G, Y_r = Y_{r-1}\cup C_r$. For each $r$, let $j_r:K_r\hookrightarrow Y_{r-1}$ denote the inclusion. Since $C_r$ is contractible and $Y_{r-1}\cap C_r = K_r$, the Mayer-Vietoris sequence for $Y_r,Y_{r-1},C_r$ yields
		$$
		H_1(K_r)\to H_1(Y_{r-1})\xrightarrow{} H_1(Y_r)\xrightarrow{} H_0(K_r)\xrightarrow{} H_0(Y_{r-1})\oplus H_0(C_r).
		$$
		Since $K_r$ is connected, the last map is injective. Exactness yields$H_1(Y_r)\cong H_1{(Y_{r-1})}/\mathrm{im}(j_r)_\ast.$
		The inclusion $j_r$ factors as 
		$$
		K_r\xhookrightarrow{\iota_r}G\hookrightarrow Y_{r-1}.
		$$
		Thus $\mathrm{im}(j_r)_\ast$ is the image in $H_1(Y_{r-1})$ of $(\iota_r)_\ast H_1(K_r)$. Iteration of this process over $r$ yields the stated formula.
	\end{proof}
	Let $X$ be a poset, $L$, $U$ and $B$ its set of minimal, maximal and nonextremal elements, respectively. For $b\in B$, set 
	$$I_{b} = \lbrace l\in L| l\leq b\rbrace,J_{b}=\lbrace u\in U| b\leq u\rbrace.$$ 
	
	We define the height of a poset as the maximal cardinality of a chain. We are now ready to give a combinatorial condition for torsion-freeness.
	\begin{theorem}\label{thm2.4}
		Let $X$ be a finite, connected poset of height $3$. If there exists a tree $T$ on $L$ such that for every $b\in B$, the induced subgraph $T[I_{b}]$ is connected, then $H_1(\Delta(X)) $ is torsion-free.
	\end{theorem}
	\begin{proof}
		Let $G$ be the bipartite graph on $L\sqcup U$ where there is an edge between two elements whenever they are comparable. Since $B$ is an antichain, we have $\Delta(X) = G\cup_{b\in B}b\ast K_{I_{b},J_{b}}$. For each $b$, let $\iota_b$ denote the inclusion map from $K_{I_{b},J_{b}}$ to $G$ and set
		$$
		S = \sum_{b\in B}(\iota_b)_\ast H_1(K_{I_{b},J_{b}})\subseteq H_1(G).
		$$
		By Lemma \ref{lem2.3} we obtain $H_1(\Delta(X))\cong H_1(G)/S$ and since $G\subseteq K_{L,U}$, the induced map
		$$H_1(G)\to H_1(K_{L,U})$$
		is injective. Thus it is enough to show that $H_1(K_{L,U})/S$ is torsion-free.
		Natural identifications $H_1(K_{L,U})\cong A_{L}\otimes A_{U}$ and $H_1(K_{I_{b},J_{b}})\cong A_{I_{b}}\otimes A_{J_{b}}$ yield $S=\sum_{b\in B}A_{I_{b}}\otimes A_{J_{b}}$ after identifying the homology groups with their images under $\iota_b$.
		Now suppose that such $T$ exists and assume an orientation on $T$. For an edge $e\in E(T)$ where $e_i$ and $e_t$ are the initial and terminal vertices, respectively, define $\alpha_e = e_t - e_i$. Notice that the elements $\alpha_e$ form a basis for $A_{L}$ and since $T[I_{b}]$ is connected, we obtain 
		$$A_{I_{b}}=\bigoplus_{e\in E(T[I_{b}])}\mathbb{Z}\alpha_e.$$
		For each $e\in E(T)$, define 
		$$C_e = \sum_{\substack{b\in B, \\e\in E(T[I_{b}])}}A_{J_b}\subseteq A_{U}.$$ Then $S = \bigoplus_{e\in E(T)}\mathbb{Z}\alpha_e\otimes C_e$, hence 
		$$(A_L\otimes A_U)/S\cong \bigoplus_{e\in E(T)}A_U/C_e.$$ By Lemma \ref{lem2.2} the group $A_U/C_e$ is torsion-free, whence the result follows.
	\end{proof}
	\section{Rectangle Covers and Side Multigraphs}
	The result that the smallest poset which admits torsion in homology has $13$ elements was originally obtained via computer verification \cite[Theorem $1.2$]{A14}. A conceptual proof which avoids extensive calculations was given in \cite{CO18}. The proof uses a method called \emph{poset splitting}. Let us briefly review the notion and relevant results as it will enable us to uncover combinatorial characterizations for the conclusiveness of a poset.

	A \textbf{beat} point in $X$ is an element which is covered by a unique element or which covers a unique element. Let $X$ be a poset on $13$ elements such that $H_1(\Delta(X))$ contains torsion. If $X$ were disconnected, then one of its connected components, necessarily having fewer than $13$ elements, would have torsion in its first homology, which is impossible. Similarly, $X$ cannot contain a beat point, because deleting such a point preserves the homotopy type and produces a poset on $12$ elements. Thus, $X$ is connected and has no beat points.
	
	Let $m_X$ and $n_X$ denote the number of maximal and minimal elements in $X$, respectively, and write $l_X$ for the number of elements that are neither maximal nor minimal. If a poset has height less than $3$, then the order complex is a graph so it cannot have torsion. Also, if $\mathrm{min}(m_X, n_X) = 2$ then $H_1(\Delta(X))$ is free. Thus it suffices to consider posets for which the height and the numbers of maximal and minimal elements are all at least $3$. By \cite[Theorem $3.9$]{CO18}, a connected finite poset satisfying \cite[Condition $(S1)$]{CO18} has free fundamental group, which implies that $H_1(\Delta(X))$ is torsion-free. Thus, $X$ does not satisfy condition $(S1)$ and \cite[Lemmas 5.3-5.5]{CO18} apply to $X$ and give 
	$$
	l_X\geq \dfrac{2m_X}{m_X-2}\text{ and }l_X\geq \dfrac{m_Xn_X}{(m_X-2)(n_X-2)} + \delta
	$$
	where $\delta =1$ if the height is greater than $3$ and $\delta = 0$ otherwise. Replacing $X$ by its dual if necessary, we may assume that $m_X\leq n_X$. Together with $m_X + l_X + n_X = 13$ and $m_X,n_X\geq 3$ the above inequalities leave only the configurations $	(3,6,4), (4,5,4),(4,4,5)	$ for the triple $(m_X,l_X,n_X)$. 
	Moreover, apart from the $(4,5,4)$ case, $B$ must be an antichain since $\delta= 1$ is incompatible with the inequalities in these cases.
	
	Although \cite[Theorem $5.10$]{CO18} is stated for finite models of $\mathbb{RP}^2$, the portions of its proof used here depend only on the fact that the poset has $13$ elements, is connected, has no beat points, and does not satisfy condition $(S1)$.
	The fact that any poset of the form $(3,6,4)$ admitting torsion is isomorphic to one of the minimal finite models of $RP^2$ is obtained in Case 1.2 of the proof of \cite[Theorem 5.10]{CO18}, whereas Case 2.2 eliminates torsion from the configuration $(4,4,5)$. Case 2.1 deals with the $(4,5,4)$ configuration and establishes that if $B$ is not an antichain, then the first integral homology group is free. However, the remainder of that case rules out a finite model of $\mathbb{RP}^2$	by means of an Euler-characteristic argument, but it does not rule out torsion in first homology. Thus, it remains to analyze the profile $(4,5,4).$
	
	The profile $(4,5,4)$ is surprisingly a lot more complicated. We will derive several combinatorial conditions that ensure torsion-freeness for such posets. We will only deal with posets without beat points as removal of such a point does not alter the homology and we already know that the poset cannot admit torsion if it has fewer than $13$ elements. 
	
	Let $B = \lbrace b_1,b_2, \dots, b_5\rbrace $. For each $i$, we set $I_i = \lbrace x\in L| x< b_i\rbrace, J_i = \lbrace x \in U| x>b_i\rbrace$. Also write $P_i = L\setminus I_i, Q_i = U\setminus J_i$ and $R_i = P_i\times Q_i$. By \cite[Lemma 5.2]{CO18}, if the missing rectangles $R_i$ do not cover the grid $L\times U$, then the hypothesis of \cite[Theorem 3.9]{CO18} is satisfied, hence the first integral homology is torsion-free. So we only have to analyze those cases where these rectangles cover $L\times U$.
	
	Suppose that $|P_i|=2$ for $i=1,2,\dots , 4$. The lower side multigraph $\Gamma_L$ is the multigraph on $L$ where there is an edge for each $P_i$. The upper side multigraph $\Gamma_U$ is defined analogously.
	The following table presents the representatives for certain types of side multigraphs that we will be using. For brevity, we omit the curly brackets and commas when denoting the sets $P_i$ and $Q_i$. Also, we use the labels $1,2,\dots$ for the maximal and minimal elements of the poset.
	\begin{table}[ht]
		\centering
		\renewcommand{\arraystretch}{1.2} 
		\setlength{\tabcolsep}{8pt}      
		
		\[
		\begin{array}{|c|c|c|}
			\hline
			\text{Notation} & \text{Representative} & \text{Description} \\
			\hline
			D& 12,12,34,34 & \text{Two disjoint doubled edges} \\
			C& 12,23,34,14 & \text{A $4$-cycle} \\
			P& 12,12,13,34 & \text{A path with one doubled edge} \\
			T& 12,13,23,14 & \text{A $3$-cycle with one extra edge attached}. \\
			\hline
		\end{array}
		\]
	\end{table}
	
	Since the poset $X$ has no beat points, for every middle element $b_i$ there are at least  two elements below and above $b_i$, respectively. Thus, $|P_i|,|Q_i|\leq 2$ and the area of each missing rectangle $R_i = P_i\times Q_i$ belongs to $\lbrace 0,1,2,4\rbrace$. The five rectangles must cover the $4\times 4$ grid $L\times U$, so the sum of their areas is at least $16$. After reordering the rectangles so that their areas are nonincreasing, we call the sequence $(|R_1|,|R_2|,\dots,|R_5|)$  the area profile of the configuration. For conciseness, we write $44422$ for $(4,4,4,2,2)$ and similarly for the other profiles. The possible area profiles are:
	$$
	44422, 44440, 44441, 44442, 44444.
	$$
	By relabeling and, if necessary, passing to the dual, one can ensure that for $i=1,2,\dots, 4$ one has $|P_i|=2$ (or $|Q_i|=2$) for each profile. So the lower (respectively upper) side multigraph can always be defined for these area profiles.
	
	We say that $X$ has a $D$-side if either $\Gamma_L$ or $\Gamma_U$ is of $D$-type (up to permutation and relabeling), and use this terminology for the other types as well. 
	Suppose $X$ is one of the profiles where there are four rectangles of area $4$ and let $R_1,R_2,\dots R_4$ denote such a quadruple. Let $\omega$ denote the number of cells that are not covered by these rectangles in the grid $L\times U$.
	\begin{theorem}\label{thm3.1}
		Let $X$ be a poset on $13$ elements of height $3$ with four maximal and four minimal elements where there are four rectangles of area $4$. Then the following holds.
		\begin{enumerate}[(i)]
			\item It is not possible to have $\omega = 1$.
			\item If $\omega = 0$, then up to duality, the pair of side multigraphs is of type $DD$ or $DC$.
		\end{enumerate}
		
	\end{theorem}
	\begin{proof}
		For the first statement, assume the contrary. Then the $4$ rectangles have one unit overlap. In the row of the uncovered cell of the grid $L\times U$, we need an overlap as each rectangle covers $0$ or $2$ cells in any row. The argument holds for columns as well, but this would imply the existence of overlaps in at least two distinct columns and rows, which makes the total overlap more than one, contradiction.
		
		If $\omega =0$, then the four rectangles partition the grid. Fix an element $x\in L$. Since there are no overlaps, exactly $2$ rectangles contribute to row $x$. Hence, $\Gamma_L$ (and analogously $\Gamma_U$) is $2$-regular. Note that a $2$-regular multigraph is either of type $D$ or $C$, so it remains to show the impossibility of the configuration $CC$. In $\Gamma_L$, let $P_i$ be an edge meeting distinct edges $P_j$ and $P_k$. The absence of overlaps forces $Q_j$ and $Q_k$ to be disjoint from $Q_i$. However, this is impossible since in a $4$-cycle for each edge there exists a unique edge that is disjoint from it. This proves the second statement.
	\end{proof}
	The tree criterion from Theorem $\ref{thm2.4}$ is particularly useful when the poset has a $D$-side.
	\begin{proposition}\label{prop3.2}
		Let $X$ be a poset on $13$ elements of height $3$ with four maximal and four minimal elements, and a $D$-side. Then $H_1(\Delta(X))$ is torsion-free.
	\end{proposition}
	\begin{proof}
		By duality, assume that $X$ has a lower $D$-side represented by $P_1=P_2=12$, $P_3=P_4=34$. Thus $I_1=I_2=34$, $I_3=I_4=12$. Start with the forest on $L$ having edges $12$ and $34$. If $I_5$ meets both $\lbrace1,2\rbrace$ and $\lbrace 3,4\rbrace$, add an edge joining an element of $I_5\cap\lbrace1,2\rbrace$ to an element of $I_5\cap\lbrace3,4\rbrace$. Otherwise add any edge joining $\lbrace 1,2\rbrace$ to $\lbrace 3,4\rbrace$. The resulting graph $T$ is a tree, and $T[I_i]$ is connected for every $i$. Theorem \ref{thm2.4} now yields the result.
	\end{proof}
	One can observe that if the area of the fifth rectangle is $0$ or $1$ then it necessarily satisfies $\omega =0$, which yields the following. 
	\begin{corollary}
		If $X$ has the area profiles $44440$ or $44441$, then $H_1(\Delta(X))$ is torsion-free.
	\end{corollary}
The $\omega=2$ case allows for more types and our next aim is to give a complete description for them. For a cell $c \in L\times U$, let $m$  denote the number of rectangles containing $c$. We call $c$ a hole if $m=0$. The excess covering at $c$ is given as $\mathrm{max}\lbrace m-1, 0\rbrace$. For an element $x\in L$, its degree $d_x$ is defined as the number of edges containing it, that is, $|\lbrace i| x\in P_i\rbrace|.$ We write $h_x$ and $e_x$ for the number of holes and the total excess covering in the row of $x$, respectively. We extend the definition to elements of $U$, analogously.
\begin{lemma}
	Let $R_1,R_2, \dots, R_4$ be four rectangles of size $4$ and suppose that $\omega =2$. Then the possible degree sequences for $\Gamma_L$ and $\Gamma_U$ are $(2,2,2,2)$ and $(3,2,2,1)$. In the latter case, both holes lie in the row or column corresponding to the degree $1$ vertex. Moreover, it is not possible to have the sequence $(3,2,2,1)$ on both sides.
\end{lemma}
\begin{proof}
	Fix $x\in L$. Notice that we have $2d_x = 4-h_x + e_x$. Since $\omega = 2$, we get $0\leq h_x,e_x\leq 2$ thus $1\leq d_x\leq 3$. The sum of degrees in a side graph is always $8$, so the possible degree sequences are $(2,2,2,2), (3,2,2,1)$ and $(3,3,1,1)$. We are to exclude the last sequence so suppose $d_x=1$. Then $h_x-e_x = 2$, which forces $h_x=2,e_x = 0$. Thus, a degree $1$ row should contain $2$ holes, meaning the sequence $(3,3,1,1)$ must have at least $4$ holes, which contradicts $\omega=2$. If both sides have this degree sequence $(3,2,2,1)$, then the contradiction follows from the fact that in addition to the row with $2$ holes, we also get a column with $2$ holes.
\end{proof}
The above degree restrictions substantially narrow the possible pairs of side graphs.
\begin{theorem}\label{thm3.4}
	Let $R_1,R_2,\dots, R_4$ be four area $4$ rectangles with $\omega = 2$. Then, up to duality, the possible pairs of side patterns are $DP,DT, PC$ and $CC$.
\end{theorem}
\begin{proof}
	By the above lemma, at least one side is $2$-regular, hence is of type $C$ or $D$. First, suppose that the lower side is $D$ and without loss of generality, let $P_1=P_2, P_3=P_4$. The rows in $P_1$ are covered only by $R_1$ and $R_2$, so the number of holes in these rows is $2(4-|Q_1\cup Q_2|)$ and in the remaining rows it similarly is $2(4-|Q_3\cup Q_4|)$. Since there are $2$ holes, one of the pairs $(Q_1,Q_2) $ and $(Q_3,Q_4)$ must be disjoint, while the other has intersection of size $1$. The disjoint pair contributes degree $1$ at every upper vertex, whereas the intersecting pair contributes $(2,1,1,0)$ so the upper side must have the sequence $(3,2,2,1)$, hence can only be $P$ or $T$, yielding the types $DP$ and $DT$.
	
	Now suppose that the lower side is $C$. From the previous argument, we know that the upper side cannot be $D$ as it forces a non-regular lower side. So we are to exclude $CT$. Write the upper $T$-side as $Q_1=12, Q_2=13, Q_3=23, Q_4 = 14$. The vertex $4$ has degree $1$, so both holes lie in this column. Column $2$ is covered by $R_1$ and $R_3$ whereas the third column is covered by $R_2$ and $R_3$, so $P_1\cup P_3=L=P_2\cup P_3$ which means that both $P_1$ and $P_2$ are disjoint from $P_3$ in the lower $4$-cycle, but such edges would have to be unique, forcing $P_1 = P_2$.
\end{proof}
It is not hard to see that the $\omega =4$ case, which can only happen in the $44444$ profile, would allow even more possibilities. The following reduction prevents this combinatorial blow up.
\begin{proposition}\label{prop3.6}
	In every configuration of profile $44444$, there exists a choice of four rectangles for which $\omega = 0$ or $\omega = 2$.
\end{proposition}
\begin{proof}
	For $i=1,2,\dots, 5$, let $\mu_i$ denote the number of cells that are uniquely covered by $R_i$. By the first assertion of Theorem \ref{thm3.1}, we already know that $\mu_i\not = 1$. If there exists an $i$ with $\mu_i =0$ or $\mu_i = 2$, the remaining quadruple of rectangles yields $\omega =0$ or $\omega = 2$, respectively. So the remaining case is $\mu_i \geq 3$ for each $i$. One can also see that $\sum_i\mu_i\leq 16$. If $\sum_i\mu_i= 16$, then every cell is uniquely covered, which is impossible with five area-$4$ rectangles. Hence $\mu_i = 3$ for each $i$. This means that $15$ cells are uniquely covered and $1$ cell belongs to all rectangles. In the row of that cell, there are $3$ elements to be covered, to which each of $5$ rectangles must contribute, but this creates more than one overlap, contradiction.
\end{proof}
\section{The Resistant Configurations}
A configuration is called \textbf{resistant} if it has no $D$-side. Proposition \ref{prop3.2} guarantees that every nonresistant configuration has torsion-free first homology. So let us first identify the resistant configurations and deal with them separately.

In the $44442$ profile, the four area-$4$ rectangles leave at most two holes and Theorem $\ref{thm3.1}$ excludes $\omega =1$. In the $44444$ profile, Proposition \ref{prop3.6} allows us to choose four rectangles for which $\omega = 0$ or $\omega =2$. It now follows from Theorems $\ref{thm3.1}$ and \ref{thm3.4} that a resistant configuration of these profiles can only be $PC$ or $CC$. In a $CC$ profile, every row vertex belongs to exactly two $P_i$. The corresponding $Q_i$ are distinct edges of the upper $4$-cycle, so their union has at least three elements. Hence every row contains at most one hole. Since there are two holes, they lie in distinct rows. By duality, they also lie in distinct columns, thus in any $CC$ configuration the two holes would have to lie in distinct rows and columns, so an area-$2$ rectangle cannot complete this grid.  Hence the only resistant configuration for $44442$ profile is $PC$, which we denote by $PC_{42}$ whereas $CC$ is possible in the $44444$ case.

Now let us deal with the $44422$ case. Observe that we do not have any overlaps. Let us first establish that the union of the two size $2$ rectangles forms a $2\times 2$ rectangle, that is, their $2$-element parts should lie in the same side. First notice that in every row and column, the number of elements that are not covered by the three size $4$ rectangles is even. It cannot be $4$ as that would force an odd number of elements to be covered in each column. So there must be exactly two rows and two columns, each containing two holes, which forces the required form for the size $2$ rectangles.

Without loss of generality, let the area $2$ rectangles be $1\times 12$ and $2\times 12$. First, suppose that one of the area-$4$ rectangles covers elements from the first two rows, so it must be $12\times 34$. The other two rectangles both have lower factor $34$, while there are two partition choices for their upper part: $\lbrace 12,34\rbrace$ or $\lbrace 13,24\rbrace$. These yield the types $1\times 12, 2\times 12, 12\times 34, 34\times 12, 34\times 34$ and $1\times 12, 2\times 12, 12\times 34, 34\times 13, 34\times 24$. Now suppose no rectangle covers elements from the first two rows, so they must occur in separate rectangles. Up to relabeling, they are $13\times 34$ and  $24\times 34$, the remaining block is then $34\times 12$. So the third type is
$1\times 12, 2\times 12, 13\times 34, 24\times 34, 34\times 12.
$ Observe that the first and third types contain a $D$-side in $\Gamma_U$. So the only resistant type of this case is the second one, which we denote by $PC_{22}$.

Unlike the others, in the $44444$ case there are several choices for the fifth rectangle to complete the grid. So let us establish these choices. Let $R_1,\dots, R_4$ be four $2\times 2$ rectangles with $\omega = 2$. A $2\times2$ rectangle $R_5$ is called a completion of this configuration if $$R_1\cup\cdots\cup R_5=L\times U.$$
\begin{prop}
	 Up to relabeling $L$ and $U$, permuting the rectangles, and taking the dual, the following holds.
	\begin{enumerate}[(i)]
		\item A configuration of type $PC$ has two completions without a $D$-side.
		\item A configuration of type $CC$ has exactly one completion.
	\end{enumerate}
\end{prop}
\begin{proof}
	Let $12\times 12, 12\times 34, 13\times 13, 34\times 24$ be a representative for the $PC$ type. The hole cells are $(4,1)$ and $(4,3)$, so the fifth rectangle $R_5 = P_5 \times Q_5$ must contain both. Consequently $Q_5 = 13$. The three possibilities for $P_5$ are $14,24$ and $34$. If $R_5 =34\times 13$, then omitting the third rectangle yields a configuration with lower $D$-side, whereas the other choices of $P_5$ do not.
	
	As observed in Proposition \ref{prop3.6}, the two holes of a $CC$ configuration must lie in distinct rows and columns. Let these holes be labeled as $(1,2)$ and $(2,1)$. We may then take the four rectangles as $13\times13, 14\times 14, 23\times24, 24\times23$. Note that the rectangles $23\times23$ and $24\times24$ would produce two more holes so the crossed pairing must be chosen. The result now follows from the fact that the holes can only be covered by $R_5=12\times 12$.
	
\end{proof}
Let $PC_{44a}$ and $PC_{44b}$ denote the two resistant $PC$ completions corresponding to $P_5 = 14$ and $P_5=24$, respectively. Also let $CC_{44}$ denote the unique resistant $CC$ completion. The following provides a useful method for eliminating torsion.
\begin{lemma}\label{lem4.2}
	Let $K$ be a finite $2$-dimensional simplicial complex. The following holds.
	\begin{enumerate}[(i)]
		\item If a $2$-simplex $\sigma = \lbrace a,b,c\rbrace$ of $K$ contains an edge $e=\lbrace a,b\rbrace$ contained in no other $2$-simplex of $K$, then $|K\setminus \lbrace\sigma, e\rbrace|$ is a strong deformation retract of $|K|$.
		\item Suppose that the $2$-simplices of $K$ can be ordered as $\sigma_1,\sigma_2,\dots,\sigma_m$ so that each $\sigma_i$ contains an edge $e_i$ belonging to none of the later simplices $\sigma_{i+1}, \sigma_{i+2},\dots ,\sigma_m$. Then $H_1(K)$ is torsion-free.
	\end{enumerate}
\end{lemma}
\begin{proof}
	Set $L = K\setminus \lbrace\sigma, e\rbrace$.  Since $K$ is $2$-dimensional and $e$ is contained in no $2$-simplex other than $\sigma$, we see that $L$ is a simplicial subcomplex of $K$. Let $\overline{\sigma}$ denote the subcomplex generated by $\sigma$. We have $K=L\cup\overline{\sigma}$ and $L\cap \overline{\sigma}=\overline{\lbrace a,c\rbrace}\cup\overline{\lbrace b,c\rbrace}$. The triangle $|\sigma|$ strongly deformation retracts onto the union $|\lbrace a,c\rbrace\cup\lbrace b,c\rbrace|$ while fixing it pointwise. Let $r_t:|\sigma|\to|\sigma|, 0\leq t\leq 1$ denote such deformation. Define $F_t:|K|\to|K|$ by
	$$
	F_t(x)=\begin{cases}
		r_t(x),& x\in|\sigma|\\
		x, &x\in |L|.
	\end{cases}
	$$
	On $|L|\cap |\sigma|$, $r_t$ is the identity, so $F_t$ is well defined and continuous. Moreover, $F_0$ is the identity map on $|K|$ and $F_1(|K|)\subseteq |L|$, which yields the first assertion.
	Let $K_0=K$ and recursively define $K_i=K_{i-1}\setminus \lbrace\sigma_i,e_i\rbrace$. At the $i$th stage, $e_i$ is contained in exactly one $2$-simplex of $K_{i-1}$. From the first part, we see that $|K_i|$ is a strong deformation retract of $|K_{i-1}|$. Thus we get $|K|\simeq |K_m|$. Since all $2$-simplices have been removed from $K_m$, we see that $H_1(K_m)$, and consequently $H_1(K)$ is torsion-free.
\end{proof}
\begin{theorem}\label{thm4.3}
	Suppose $X$ belongs to one of the configurations $PC_{42}, PC_{44a}, PC_{44b}, CC_{44}$. Then $H_1(\Delta(X))$ is torsion-free.
\end{theorem}
\begin{proof}
	Let $T_i=\lbrace [x,b_i,z]| x\in I_i, z\in J_i\rbrace$ denote the set of $2$-simplices through $b_i$. Let us first deal with the $CC_{44}$ case. The five $I_i\times J_i$ rectangles are $24\times24, 23\times 23, 14\times 13, 13\times14$ and $34\times 34$. Notice that the cells $(2,4),(2,3),(1,3),(1,4)$ belong to exactly one of the five rectangles. So for $i= 1,2,\dots, 4$, there exists an element $(x_i,z_i)\in L\times U$ uniquely contained in the $i$th rectangle. Consequently, the edge $[x_i,z_i]$ occurs in no $2$-simplex outside $T_i$.  For each $i$, let $[x_i,b_i,z_i], [x_i,b_i,t_i],[y_i,b_i,z_i],[y_i,b_i,t_i]$ be an ordering of $T_i$. The edges $[x_i,z_i], [x_i,b_i], [b_i,z_i]$ and $[y_i,b_i]$ from each $2$-simplex occur in no later one. Choose any ordering on $T_5$. The distinguishing edge for each $2$-simplex $[x,b_i,z]\in T_5$ can simply be chosen as $[x,z]$. The resulting concatenation of the families $T_i$ satisfies the hypothesis of the above lemma.
	
	For $PC$ types, the first four $I_i\times J_i$ rectangles are $34\times 34, 34\times 12, 24\times 24, 12\times 13$, whereas the fifth for each case is given by
	$$
	I_5\times J_5 =\begin{cases}
		123\times 24, &X=PC_{42}\\
		23\times 24 &X=PC_{44a}\\
		13\times 24 &X=PC_{44b}.
	\end{cases}
	$$
	
	Notice that the cells $(3,3)$ and $(3,1)$ belong to exactly one of the five rectangles in each configuration. Using the internal ordering for $T_i$'s given for the previous case, the required simplex-edge ordering for $T_1$ and $T_2$ can be obtained. We are left with $T_3,T_4$ and $T_5$. Notice that the cells $(4,2)$ and $(1,1)$ belong to the third and fourth rectangles, respectively, but do not occur in a later one. Thus, we can concatenate $T_3$ and $T_4$ with the above internal ordering to the existing sequence. Similar to the above, in any ordering for $T_5$, we can find the distinguishing edges which occur in no later simplex. Appending $T_5$ completes the required ordering.
	
\end{proof}
In order to complete the proof of Theorem $\ref{thm1.1}$, it remains to eliminate torsion from the $PC_{22}$ configuration, for which the $I_i\times J_i$ rectangles are given by:
$$
34\times 12, 12\times24, 12\times 13, 234\times 34, 134\times 34.
$$
Note that it is not possible to order every $2$-simplex of $PC_{22}$ in a way that satisfies the hypothesis of the second part of Lemma \ref{lem4.2}. The problem is with the last two rectangles. However, it is still possible to repeatedly apply the first part of the same lemma and obtain a simplified structure. In order to deal with this remaining complex, we shall prepare some topological machinery. For a simplicial complex $K$ and a vertex $v\not\in K$, let $v\ast K$ denote the cone on $K$ with apex $v$.
\begin{lemma}
	For a nonempty, finite, connected simplicial complex $K$, the following holds.
	\begin{enumerate}[(i)]
		\item For vertices $v$ and $w$ not belonging to $K$, set $S=(v\ast K)\cup (w\ast K)$. Then $H_1(S)=0$.
		\item If $K$ is $2$-dimensional and has a subcomplex $S$ that contains every $2$-simplex of $K$ and satisfies $H_1(S)=0$, then $H_1(K)$ is torsion-free.
	\end{enumerate}
\end{lemma}
\begin{proof}
	Set $A=v\ast K, B=w\ast K$. We have $H_1(A)=H_1(B)=0$ and $H_0(A)\cong H_0(B)\cong \mathbb{Z}$ and since $K$ is connected, we also have $H_0(K)\cong \mathbb{Z}$. Notice that $A\cap B = K$. The Mayer-Vietoris sequence for $S,A$ and $B$ contains
	$$
	H_1(A)\oplus H_1(B)\to H_1(S)\to H_0(K)\to H_0(A)\oplus H_0(B).
	$$
	 Since $K$, $A$ and $B$ are connected, both inclusion maps induce isomorphisms on $H_0$, hence the last map is injective. Exactness now yields the first assertion.
	 
	For the second statement, consider the long exact sequence for the pair $(K,S)$ which contains
	$$
	H_1(S)\to H_1(K)\to H_1(K,S)\to H_0(S)\to H_0(K).
	$$
	Since $S$ contains every $2$-simplex of $K$, we have  $C_2(K,S)=0$, thus $H_1(K,S)$ is a subgroup of $C_1(K,S)$, and is free. Because $H_1(S)=0$, the exactness yields that the natural map $H_1(K)\to H_1(K,S)$ is an injection, from which the result follows.
\end{proof}
\begin{theorem}
	Suppose $X$ is a poset of type $PC_{22}$. Then $H_1(\Delta(X))$ is torsion-free.
\end{theorem} 
\begin{proof}
	Retain the notation $T_i$ from Theorem \ref{thm4.3}. In $T_1$, pair the four simplices with the edges $$(3,1), (3,2), (4,1), (4,2).$$ In $T_2$, first use the edges $(1,2)$ and $(2,2)$, and then use the edges $(1,b_2)$ and $(b_2,4)$. Similarly, in $T_3$, first use $(1,1)$ and $(2,1)$, then $(1,b_3)$ and $(b_3, 3)$. Consequently, removal of these families will preserve the homotopy type, so we are left with $T_4$ and $T_5$. Notice that in $T_4$, the edges $(2,3)$ and $(2,4)$ can be paired with the simplices $[2,b_4,3]$ and $[2,b_4,4]$, respectively and in $T_5$, $(1,3)$ and $(1,4)$ can be paired with $[1,b_5,3]$ and $[1,b_5,4]$, respectively. We can further remove these simplices, which leaves us with the complex $Y$, whose $2$-simplices $[x,b_i,z]$ are given by $x,z\in\lbrace 3,4\rbrace,i\in\lbrace 4,5\rbrace$. 
	
	Notice that the remaining edges between the minimal and maximal elements form the bipartite graph $K_{2,2}$, hence all $2$-simplices of $Y$ form the subcomplex $S=(b_4\ast K_{2,2})\cup (b_5\ast K_{2,2})$. Since $K_{2,2}$ is connected, the first part of the above lemma yields $H_1(S)=0$. Moreover, $S$ contains every $2$-simplex of $Y$. Thus the second part of the above lemma shows that $H_1(Y)$ is torsion-free. Since $Y$ is obtained from $\Delta(X)$ by successively removing the simplex-edge pairs satisfying Lemma \ref{lem4.2}(i), it follows that $H_1(\Delta(X))$ is torsion-free.
\end{proof}

\begin{proof}[Proof of Theorem \ref{thm1.1}]
	We can now complete the proof of Theorem $1.1$. Let $X$ be a poset on $13$ elements whose first integral homology contains torsion. By the preceding reductions, up to duality, $X$ has one of the profiles $(3,6,4), (4,5,4)$, or $(4,4,5)$. The first profile gives one of the two minimal finite models of $\mathbb{RP}^2$, while the third profile is torsion-free. In the remaining profile $(4,5,4)$, the nonextremal elements form an antichain. Sections 3 and 4 show that every possible rectangle configuration of this profile has torsion-free first homology. Therefore the only $13$-element posets with torsion in first homology are the two dual minimal finite models of $\mathbb{RP}^2$. Conversely, both of these posets have first homology isomorphic to $\mathbb{Z}/2\mathbb{Z}$. This proves Theorem $1.1$.
\end{proof}
Theorem \ref{thm1.1} and the above developed methods raise the natural question of whether for each $n\geq 13$, there exists a poset $P$ on $n$ elements such that $H_1(\Delta(P))$ contains torsion whereas for each $x\in P$, the group $H_1(\Delta(P\setminus \lbrace x\rbrace))$ is torsion-free. We do not know whether such posets exists for every $n\geq 14$.

\section{Further Criteria for Outer Derivations}
Having completed the homological analysis, we return to transitive functions and derive two further criteria for the existence of outer derivations.
Let $P$ be a finite poset. A directed path in the Hasse diagram of $P$ is a sequence $x_0\prec x_1\prec\cdots\prec x_r$ where each relation is a cover relation.
 Two distinct directed paths with the same initial and terminal vertices form a pair of parallel paths. Each such pair produces a consistency relation by equating the sums of the edge values along the two paths.
 
Let $\mathcal{M}_P$ denote the matrix whose columns are indexed by the edges of the Hasse diagram of $P$ and the rows are the coefficient vectors of the consistency relations associated with pairs of parallel paths. Write $C(P)$ for the number of connected components of $P$ and let $E(P)$ and $V(P)$ denote the edge and vertex sets of the Hasse diagram of $P$, respectively. Also let $\mathrm{Pot}(P,k)$ and $\mathrm{Tr}(P,k)$ denote the $k$-vector spaces of potential and transitive functions on $P$, respectively. Notice that $\mathrm{Pot}(P,k)\subseteq\mathrm{Tr}(P,k)$. The following gives a complete criterion for the absence of outer derivations.

\begin{theorem}
Let $P$ be a finite poset, $k$ a field. The following are equivalent:
\begin{enumerate}[(i)]
\item $\mathrm{Pot}(P,k)=\mathrm{Tr}(P,k)$.
\item $|E(P)| - \mathrm{rank}_k(\mathcal{M}_P) = |V(P)| -C(P).$
\end{enumerate}
\end{theorem}
\begin{proof}
Every transitive function is uniquely determined by its values on the cover relations. Conversely, an assignment of values to the edges of the Hasse diagram extends to a transitive function if and only if these sums are independent of the chosen path, or equivalently, if the sums agree along every pair of parallel paths. Since these equalities are encoded by the rows of $\mathcal{M}_P$, restriction to the edges gives an isomorphism $$\mathrm{Tr}(P,k)\cong\mathrm{ker}(\mathcal{M}_P),$$ consequently $$\mathrm{dim}_k(\mathrm{Tr}(P,k))=|E(P)| - \mathrm{rank}_k(\mathcal{M}_P).$$
It is enough to show that when $P$ is connected, one gets $\mathrm{dim}_k(\mathrm{Pot}(P,k)) = |V(P)|-1$. 
Let $N$ be the $k$-vector space of all functions on $V(P)$. The $k$-dimension of $N$ is clearly $|V(P)|$. Any element $\varphi\in N$ yields a potential function $\delta_\varphi$ where 
    $$
    \delta_\varphi(x,y) = \varphi(y) - \varphi(x)
    $$
for all $x,y\in P, x< y$. The map from $N$ to $\mathrm{Pot}(P,k)$ defined by $\varphi\mapsto \delta_\varphi$ is a homomorphism of the $k$-modules and is surjective. Since $P$ is connected, the kernel consists of constant functions, and therefore has dimension $1$, whence the formula in the connected case follows. The statement for $C(P)>1$ can be obtained by applying this argument to each component.
\end{proof}
Every graph $G$ with $V$ vertices, $E$ edges and $C$ connected components satisfies $E \geq V-C$ and equality holds if and only if $G$ is acyclic. This yields the following known result (\cite[Corollary 3.5]{FP21}), which originally was obtained via a different approach.
\begin{corollary}\label{cor5.2}
Let $P$ be a finite poset with no parallel paths. Then $P$ admits no outer derivations over a field $k$ if and only if the underlying undirected graph of its Hasse diagram is a forest.
\end{corollary}

A \textbf{crown} is a poset $C_n = \lbrace x_1,x_2,\dots,x_n, y_1,y_2,\dots y_n\rbrace$ where $n\geq 2$ and for $1\leq i\leq n-1$, we have $x_i< y_i, x_i < y_{i+1}$ and $x_n <y_n, x_n <y_1$. The crown $C_n$ can be visualized as follows.
\begin{figure}[ht] 
    \centering
\begin{tikzpicture}[scale = 0.8,
    x=1.5cm, y=1.5cm,           
    dot/.style={shape=circle, fill=black, inner sep=0pt, minimum size=4pt}, 
    edge/.style={thick}         
]

    \foreach \i in {1, 2, 3} {
        \node[dot] (x\i) at (\i-1, 0) {}; \node[below=3pt] at (x\i) {$x_{\i}$};
        \node[dot] (y\i) at (\i-1, 1) {}; \node[above=3pt] at (y\i) {$y_{\i}$};
    }
	\node[dot] (y4) at (3, 1) {}; \node[above=3pt] at (y4) {$y_4$};
    \node (dots_x) at (3, 0) {$\dots$};
    \node (dots_y) at (4, 1) {$\dots$};

    \node[dot] (xn) at (5, 0) {}; \node[below=3pt] at (xn) {$x_n$};
    \node[dot] (yn) at (5, 1) {}; \node[above=3pt] at (yn) {$y_n$};
	\node[dot] (xn1) at (4, 0) {}; \node[below=3pt] at (xn1) {$x_{n-1}$};

    \draw[edge] (x1) -- (y1);
    \draw[edge] (x2) -- (y2);
    \draw[edge] (x3) -- (y3);
    \draw[edge] (x3) -- (y4);
    \draw[edge] (xn) -- (yn);
	\draw[edge] (xn1) -- (yn);
    \draw[edge] (x1) -- (y2);
    \draw[edge] (x2) -- (y3);
     \draw[edge] (xn) -- (y1);

\end{tikzpicture}
\caption{The crown $C_n$ on $2n$ elements.}
\end{figure}

A cycle in a poset $P$ is a sequence $(x_1,\cdots,x_n)$ of elements where $x_i\in P$ such that $x_n=x_1$ and for all $1\leq i\leq n-1$, $x_i$ and $x_{i+1}$ are comparable.
Given a transitive function $f$, define
$$
\tilde{f}(x,y) = \begin{cases}
			f(x,y), & \text{if }x<y\\
            -f(y,x), & \text{if } y<x.
		 \end{cases}
$$
For a cycle $C = (x_1,x_2,\cdots, x_n)$ the circulation $v_f(C)$ of $f$ on $C$ is defined as
$$
v_f(C) = \sum_{i=1}^{n-1}\tilde{f}(x_i,x_{i+1}).
$$
It is known that containment of a crown is a requirement for the non-vanishing of the first cohomology. In the context of transitive systems, this fact is given as \cite[Theorem $1.12$]{CJW22}. The following is a strengthening.
\begin{theorem}\label{thm5.3}
    Suppose that $P$ is a finite poset admitting a nonpotential transitive function $f$ with codomain $G$. Then $P$ contains a crown $C$ satisfying 
    $$
    v_f(C)\not=0
    $$
    such that the elements of $C$ have neither a common upper bound nor a common lower bound in $P$.
\end{theorem}
\begin{proof}
By \cite[Theorem $1.5$]{CJW22}, the existence of a nonpotential transitive function implies the existence of a cycle with nonzero circulation. Choose such a cycle $C$ of minimal length. By minimality, no vertex is visited more than once (except the initial one) in $C$ as in that case, we can divide $C$ into two shorter cycles and obtain a smaller cycle with nonzero circulation. Similarly, no two nonconsecutive vertices $x_i,x_j$ are comparable. Indeed, such a comparability would give a chord dividing $C$ into two shorter cycles. The chord occurs with opposite signs in their circulations, so the circulations of the two cycles sum to $v_f(C)$. Consequently, one of them would have nonzero circulation, contradicting the minimality of $C$. Thus $C$ is chordless, in particular, no three consecutive vertices of $C$ can form a chain,  hence $C$ is a crown. 

Let $C = (x_1,x_2,\cdots, x_{2k + 1})$ where $x_1 = x_{2k+1}$. Now suppose that the elements of $C$ have a common upper bound $M$ in $P$. Then for $i =1,2,\cdots, 2k$ we have 
$$
\tilde{f}(x_i,x_{i+1}) = \tilde{f}(x_{i},M) +\tilde{f}(M,x_{i+1})= f(x_{i},M) -f(x_{i+1},M).
$$
 This yields
    \begin{align*}
        v_f(C) &=  \sum_{i=1}^{2k}\tilde{f}(x_i,x_{i+1})\\
        &= f(x_1,M)-f(x_2,M)+f(x_2,M)-f(x_3,M) +\cdots + f(x_{2k},M) - f(x_{2k+1},M) \\
        &= f(x_1,M) - f(x_{2k+1},M)=0.
    \end{align*}
So the circulation must vanish, which is a contradiction. Hence, $C$ has no common upper bound. The nonexistence of a lower bound follows from duality.
\end{proof}

\section*{Acknowledgements}
This study was supported by the Scientific and Technological Research Council of Türkiye (TUBITAK) under Grant Number 125F395. The authors thank TUBITAK for its support.
\section*{Data Availability}
No data were used for the research described in the article.
\section*{Conflict of Interests}
The authors declare that they have no conflict of interests.

\end{document}